\documentclass[preprints,article,accept,pdftex,moreauthors]{my-mdpi} 

\firstpage{1} 
\pubvolume{11}
\issuenum{8}
\articlenumber{379}
\pubyear{2022}
\copyrightyear{2022}
\externaleditor{Academic Editor: Valery Y. Glizer} 
\datereceived{6 July 2022} 
\daterevised{22 July 2022} 
\dateaccepted{29 July 2022} 
\datepublished{31 July 2022} 
\hreflink{https://\\doi.org/10.3390/axioms11080379} 
\doinum{10.3390/axioms11080379}

\makeatletter      
\renewcommand{\maketag@@@}[1]{\hbox{\m@th\normalsize\normalfont#1}}%
\makeatother
\pdfoutput=1

\Title{Minimum Energy Problem in the Sense of Caputo for Fractional Neutral Evolution Systems in Banach Spaces}

\TitleCitation{Minimum Energy Problem in the Sense of Caputo for Fractional Neutral Evolution Systems in Banach Spaces}

\Author{Zoubida Ech-chaffani $^{1,\dagger}$\orcidA{}, 
Ahmed Aberqi $^{2,\dagger}$\orcidB{}, 
Touria Karite $^{3,}$*$^{,\dagger}$\orcidC{} 
and Delfim F. M. Torres $^{4,\dagger}$\orcidD{}}

\AuthorNames{Zoubida Ech-chaffani, Ahmed Aberqi, Touria Karite and Delfim F. M. Torres}

\AuthorCitation{Ech-chaffani, Z.; Aberqi, A.; Karite, T.; Torres, D.F.M.}

\address{$^{1}$ \quad LAMA Laboratory, Faculty of Sciences Dhar El Mahraz, 
Sidi Mohamed Ben Abdellah University, \linebreak  
Fez--Atlas BP 1796, Morocco; echzoubida@gmail.com\\
$^{2}$ \quad LAMA Laboratory, Department of Electrical Engineering \& Computer Science, 
National School of Applied Sciences, Sidi Mohamed Ben Abdellah University, 
Avenue My Abdallah Km 5 Route d'Imouzzer, \linebreak  Fez BP 72, Morocco; ahmed.aberqi@usmba.ac.ma\\
$^{3}$ \quad Laboratory of Engineering, Systems and Applications, 
Department of Electrical Engineering \& Computer Science, 
National School of Applied Sciences, Sidi Mohamed Ben Abdellah University, 
Avenue My Abdallah Km 5 Route d'Imouzzer, Fes BP 72, Morocco\\
$^{4}$ \quad Center for Research and Development in Mathematics and Applications (CIDMA), 
Department of Mathematics, University of Aveiro, 3810-193 Aveiro, Portugal; delfim@ua.pt}

\corres{Correspondence: touria.karite@usmba.ac.ma}

\firstnote{These authors contributed equally to this work.}

\abstract{We investigate a class of fractional neutral evolution 
equations on Banach spaces involving Caputo derivatives. 
Main results establish conditions for the controllability 
of the fractional-order system and conditions for existence
of a solution to an optimal control problem of
minimum energy. The results are proved with 
the help of fixed-point and semigroup theories.}



\begin{document}
	

\section{Introduction}
\label{sec1}

A neutral system is a system where time-delays play an important role. 
Precisely, such delays appear in both state variables and their derivatives.
A delay in the derivative is called ``neutral'', which makes the system 
more complex than a classical one where the delays only occur in the state. 
Neutral delays do not only occur in physical systems, but~they also 
appear in control systems, where they are sometimes added to improve 
the performance. For~instance, a~wide range of neutral-type control
systems are expressed by
\begin{equation}
\label{eq:introd:ns1}
\frac{d}{dt}\left[y(t)-Ky_{t}\right] = Ly_{t} + Bu(t), 
\quad t\geq 0,\quad y_{0}(\cdot) = f_{0}(\cdot),
\end{equation}
where $y_{t} : [-1,0] \rightarrow \mathbb{C}^{n}$ 
is defined by $y_{t}(s) = y(t + s)$;
for $f \in H^{1}([-1,0],\mathbb{C}^{n})$,
the difference operator $K$ is given by
$Kf = A_{-1}f(-1)$ with $A_{-1}$ a constant $n\times n$ matrix.
The delay operator $L$ is defined by
$$
Lf = \displaystyle\int_{-1}^{0}\left[A_{2}(\theta)f'(\theta) 
+ A_{3}(\theta)f(\theta)\right]d\theta
$$
with $A_{2}$ and $A_{3}$ $n\times n$ matrices whose elements belong 
to $L_{2}(-1,0)$; $B$ is a constant $n\times r$ matrix;
and the control $u$ is an $L_{2}$-function~\cite{Rabah14}.

Nowadays, many researchers have investigated neutral differential equations 
in Banach spaces~\cite{Chandrasekaran15,Huang22,Harisa22}. This interest 
is explained by the fact that neutral-argument differential equations have 
interesting applications in real-life problems: they appear, e.g.,~while 
modeling networks containing lossless transmission lines or in super-computers. 
Moreover, second-order neutral equations play an important role in 
automatic control and in aeromechanical systems, where inertia plays 
a central role~\cite{Baker2008,Baculikova17,Bazighifan19}.

Controllability plays an inherent crucial role in finite and infinite-dimensional 
systems, being one of the primary concepts in control theory, along with observability 
and stability. This concept has also attracted many authors; 
see, for~instance,~\cite{Book:Lions:1971,K9,Karite21}.

In the last two decades, several researchers have been interested 
in exploring the concept of controllability for fractional 
systems~\cite{Ge18,MR3895586,Almeida19}. This is natural because 
fractional differential equations are considered a valuable 
tool in modeling various real-world dynamic 
systems, including~physics, biology, socio-economy, chemistry 
and engineering~\cite{Hilfer2000,Magin2006,Tabanfar22}.

It turns out that system \eqref{eq:introd:ns1} can also
be studied in the fractional sense, e.g.,~being expressed by
$$
\left\{\begin{array}{ll}
{ }^{C} D_{t}^{q}\left[y(t)-Ky_{t}\right] 
=Ly(t)+Bu(t), & t \in [0, T],\\
y_{0}(\cdot)=f_{0}(\cdot),
\end{array}
\right. 
$$
where ${ }^{C} D_{t}^{q}$ denotes the Caputo fractional derivative of order $q$.
The existence of solutions to fractional differential equations 
for neutral systems involving Caputo or other fractional operators,
\textls[-15]{like Riemann--Liouville fractional derivatives, 
has been paid much attention~\cite{24,25,26}.} Recently, 
some achievements regarding the existence and uniqueness of mild solutions 
to fractional stochastic neutral differential systems
in a finite dimensional space have been made~\cite{Ahmadova20}. 
Other works are consecrated to demonstrate existence 
of a mild solution for neutral fractional inclusions 
of the Sobolev type~\cite{Kavitha21}.

In~\cite{28}, Sakthivel~et~al. examined the exact controllability 
of fractional differential neutral systems by establishing 
sufficient conditions via a fixed-point analysis approach. 
Later on, Sakthivel~et~al. investigated the weak controllability 
of fractional dynamical systems of order $1 < q < 2$ using 
sectorial operators and Krasnoselskii's fixed-point theorem~\cite{Sakthivel13}. 
Using the same techniques as the previous authors, Qin~et~al. have studied
the controllability and~optimal control of~fractional dynamical systems 
of order $1 < q < 2$ in Banach spaces~\cite{Qin15}. Yan and Jia used 
stochastic analysis theory and fixed-point theorems 
with the strongly continuous $\alpha$-order cosine family 
to study an optimal control problem for a class of stochastic 
fractional equations of order $\alpha\in (1,2]$ in Hilbert spaces~\cite{Yan18}. 
In 2021, Zhou and He obtained, via the contraction principle 
and Shauder's fixed-point theorem, a~set of sufficient conditions 
for the exact controllability of a class of fractional systems~\cite{31}. 
More recently, Xi~et~al. studied the approximate controllability of fractional 
neutral hyperbolic systems using Sadovskii's fixed point theorem while constructing 
a Cauchy sequence and a control function~\cite{Xi22}. Dineshkumar~et~al. addressed 
the problem of approximate controllability for neutral stochastic fractional 
systems in the sense of Hilfer, treating the problem using Schauder's fixed-point 
theorem and extending the obtained results to the case of nonlocal 
conditions~\cite{Dineshkumar21}. In~\cite{Ma22}, Ma~et~al. analyzed the weak controllability 
of a fractional neutral differential inclusion of the Hilfer type in Hilbert spaces 
using Bohnenblust--Karlin's fixed point theorem. The~concept of complete controllability 
is studied in~\cite{Wen22} by Wen and Xi, where they establish sufficient conditions 
to assure this type of controllability. 

Here, we let $(X, \mid\cdot\mid)$ be a Banach space, 
and we denote the Banach space of continuous functions 
by $ \mathcal{C}(0,T; X)$ with the norm 
$\vert x\vert =\displaystyle\sup_{t \in J}\vert x(t)\vert$.
Our main goal is to explore the concepts of controllability 
and optimal control for the following general evolution 
fractional system:
\begin{equation}
\label{system1}
\left\{\begin{array}{ll}
{ }^{C} D_{t}^{\nu}\left[x(t)-h\left(t, x_{t}\right)\right] 
=\mathcal{A} x(t)+\mathcal{B}u(t), & t \in (0, T],\\
x(0)=x_{0} \in D(\mathcal{A}),
\end{array}
\right. 
\end{equation}
where ${ }^{C} D_{t}^{\nu}$ denotes the fractional derivative 
of order $\nu\in(0,1)$ in the sense of Caputo, 
$h:[0,T]\times \mathcal{C}(0, T ; X)\rightarrow X$ 
is a given continuous function,
and the dynamic of the system 
$\mathcal{A}:D(\mathcal{A})\subseteq X\rightarrow X$ is a linear, 
closed operator with dense domain $D(\mathcal{A})$ generating 
a compact and uniformly bounded $C_{0}$ semigroup 
$\{\mathcal{T}(t)\}_{t \geq 0}$ on $X$. The~control function 
$u(\cdot)$ is given in $L^{2}(0,T;U)$, with~$U$ a reflexive Banach space,
and the control operator $\mathcal{B}\in \mathcal{L}(U,X)$ is a linear 
continuous bounded operator, i.e.,~there exists a constant $M_{1}>0$ 
such that\vspace{-6pt}
\begin{equation}
\label{eq2}
\displaystyle \vert \mathcal{B}\vert \leq \mathit{M}_{1}.
\end{equation}
Our main aim is to be able to obtain a set of sufficient 
conditions assuring the controllability of system \eqref{system1} 
and, afterwards, to~consider an associated optimal control problem 
and prove existence of a solution.

The rest of this paper is organized as follows. In~Section~\ref{sec2}, 
the definitions of Caputo fractional derivative and mild solutions for 
system \eqref{system1} are recalled. Our main result on the controllability 
of \eqref{system1} is proved in Section~\ref{sec3}.
In Section~\ref{sec4}, we prove the existence of a control 
giving minimum energy on a closed convex set of admissible controls. 
Section~\ref{sec5} is consecrated to the analysis of a concrete example, 
illustrating the applicability of our main results. 
We end with Section~\ref{sec6}, which contains conclusions
and points out some possible future directions of research.


\section{Background}
\label{sec2}

\textls[-15]{In this section, basic definitions, notations, and~lemmas are introduced 
to be used throughout the paper. In~particular, we recall the main 
properties of fractional calculus~\cite{I.P,A.A}} 
and useful properties of semigroup theory~\cite{Pazy}.

Throughout the paper, let $\displaystyle \mathcal{A}$ be 
the infinitesimal generator of a compact and uniformly bounded 
$\displaystyle C_{0}$ semi-group $\displaystyle \{\mathcal{T}(t)\}_{t \geq 0}$. 
Let $\displaystyle 0\in\varrho(\mathcal{A})$, where $\displaystyle \varrho(\mathcal{A})$ 
denotes the resolvent of $\mathcal{A}$. Then, for~$\displaystyle 0\leq\mu\leq 1$, 
the fractional power $\mathcal{A}^{\mu}$ is defined as a closed linear 
operator on its domain $\displaystyle D(\mathcal{A}^{\mu})$.
For a compact semi-group $\displaystyle \{\mathcal{T}(t)\}_{t \geq 0}$, 
the following properties are useful in this~paper:
\begin{enumerate}
\item[(i)] There exists $\mathit{M}_{T}\geq1$ such that
\begin{equation}
\label{M_{T}}
\displaystyle \mathit{M}_{T}=\sup_{t\geq 0}\vert \mathcal{T}(t)\vert ;
\end{equation}
\item[(ii)] For any $\displaystyle \mu \in (0,1]$, 
there exists $\displaystyle \mathbb{L}_{\mu}>0$ such that
\begin{equation}
\label{C}
\displaystyle \vert \mathcal{A}^{\mu}\mathcal{T}(t)\vert 
\leq\displaystyle\frac{\mathbb{L}_{\mu}}{t^{\mu}},
\quad 0\leq t\leq T.
\end{equation} 
\end{enumerate}

Now we recall the notion of a Caputo fractional~derivative.

\begin{Definition}[See~\cite{A.A}]
The left-sided Caputo fractional derivative of order 
$\displaystyle \nu>0$ of a function $z\in L^{1}([0,T])$ is
\begin{equation}
\label{D3}
\displaystyle { }_{0}^{C} \mathit{D}_{t}^{\nu} z(t)
=\frac{1}{\Gamma(\kappa-\nu)} \int_{0}^{t}(t-s)^{\kappa-\nu-1} 
\frac{d^{\kappa}}{d s^{\kappa}} z(s) d s,
\end{equation}
where $t\geq 0$, $\kappa-1<\nu<\kappa$, $\kappa\in\mathbb{N}$, 
and  $\displaystyle \Gamma(\cdot)$ is the gamma function.
\end{Definition}

Using the probability density function and its Laplace transform~\cite{Y.F} 
(see also~\cite{M.M, Kumar}), we recall the definition of a mild solution 
for system \eqref{system1}.

\begin{Definition}[See~\cite{Y.F}]
Let $\displaystyle u\in U$ for $\displaystyle t\in] 0, T]$. 
A function $\displaystyle x \in \mathcal{C}(0, T ; X)$ is said 
to be a mild solution of system $ \eqref{system1} $ if~\begin{equation}
\label{mild solution}
\begin{aligned} 
x(t,u) &=\mathit{S}_{\nu}(t)\left[x_{0}-h\left(0, x_{0}\right)\right]
+h\left(t, x_{t}\right)+\int_{0}^{t}(t-s)^{\nu-1} \mathcal{A} 
\mathit{K}_{\nu}(t-s) h\left(s, x_{s}\right) \mathrm{d} s \\ 
&\quad +\int_{0}^{t}(t-s)^{\nu-1} \mathit{K}_{\nu}(t-s) 
\mathcal{B}u(s)\mathrm{d} s,
\end{aligned}
\end{equation}
where $\displaystyle \mathit{S}_{\nu}(\cdot)$ and  
$\displaystyle \mathit{K}_{\nu}(\cdot)$ are the 
characteristic solution operators defined by
$$
\mathit{S}_{\nu}(t) 
=\int_{0}^{\infty} \phi_{\nu}(\Theta) 
\mathcal{T}\left(t^{\nu} \Theta\right)  \mathrm{~d} \Theta 
\quad \mbox{and} \quad 
\mathit{K}_{\nu}(t) 
=\nu \int_{0}^{\infty} \Theta \phi_{\nu}(\Theta) 
\mathcal{T}\left(t^{\nu} \Theta\right)  \mathrm{~d}\Theta
$$
with 
$$
\displaystyle \phi_{\nu}(\Theta)
=\frac{1}{\nu}\Theta^{-1-\frac{1}{\nu}}
\psi_{\nu}\left(\Theta^{-\frac{1}{\nu}}\right)
$$
and
$$
\psi_{\nu}(\Theta)=\frac{1}{\pi} \sum_{n=1}^{\infty}(-1)^{n-1} \Theta^{-{\nu n-1}} 
\frac{\Gamma(n \nu+1)}{n!} \sin (n\pi \nu), 
\quad \Theta \in(0, \infty),
$$  
the probability density. In~addition, we have 
$$
\int_{0}^{\infty} \psi_{\nu}(\Theta) d \Theta=1 
\text{ and } \int_{0}^{\infty} \Theta^{\Lambda} 
\phi_{\nu}(\Theta) d \Theta=\frac{\Gamma(1+\Lambda)}{\Gamma(1+\nu \Lambda)}, 
\quad \Lambda \in[0,1].
$$
\end{Definition}

\begin{Remark}
The solution $x(t,u)$ of \eqref{system1} is considered in the weak sense, 
and, when there are no ambiguities, it is denoted by $x_{u}(t)$. We denote 
by $x_{u}(T)$ the mild solution of system \eqref{system1} at the final time $T$.
\end{Remark}

The following properties of $\displaystyle \mathit{S}_{\nu}(\cdot)$ 
and  $\mathit{K}_{\nu}(\cdot)$ will be used throughout the~paper.

\begin{Lemma}[See~\cite{Y.F}] \ \ 
\label{YF}
\begin{enumerate}

\item For any $\displaystyle t\geq 0$, the~operators 
$\displaystyle \mathit{S}_{\nu}(t)$ 
and $\displaystyle \mathit{K}_{\nu}(t)$ are linear and bounded, i.e.,~
$$
\left\vert \mathit{S}_{\nu}(t) y\right\vert 
\leq \mathit{M}_{T}\vert y\vert 
\quad \text{ and } \quad 
\left\vert \mathit{K}_{\nu}(t) y\right\vert 
\leq \frac{\nu \mathit{M}_{T}}{\Gamma(1+\nu)}\vert y\vert
$$
for any $\displaystyle y\in X$
where $\displaystyle \mathit{M}_{T}=\sup_{t\geq 0}\vert \mathcal{T}(t)\vert$.

\item For $\displaystyle t>0$, if~$\displaystyle \mathcal{T}(t)$ is compact, 
then $\displaystyle \mathit{S}_{\nu}(t)$ and $\displaystyle \mathit{K}_{\nu}(t)$ 
are both compact operators.
\end{enumerate}
\end{Lemma}

\begin{Lemma}[See~\cite{Y.F}]
\label{2.2}
For any $\displaystyle x\in X$, $\displaystyle \varsigma\in (0,1)$
and $\displaystyle \mu\in (0,1]$ we~have
\begin{itemize}

\item[(i)] $\displaystyle \mathcal{A} \mathit{K}_{\nu}(t) x
=\mathcal{A}^{1-\varsigma} \mathit{K}_{\nu}(t) \mathcal{A}^{\varsigma} x$, 
$0 \leq t \leq a$;

\item[(ii)]$\displaystyle \left\vert \mathcal{A}^{\mu} 
\mathit{K}_{\nu}(t)\right\vert  
\leq \displaystyle\frac{\nu \mathbb{L}_{\mu}}{t^{\nu \mu}} 
\frac{\Gamma(2-\mu)}{\Gamma(1+\nu(1-\mu))}$, $0<t \leq a$.
\end{itemize}
\end{Lemma}


\section{Controllability}
\label{sec3}

Following~\cite{F.Ge}, let us define the meaning 
of controllability for our system \eqref{system1}.

\begin{Definition}
System \eqref{system1} is said to be controllable in $X$ on $[0,T]$ 
if for any given initial state $x_{0}\in X$ and any desired final 
state $x_{d}\in X$, there exists a control $u(\cdot)\in L^{2}(0,T;U)$ 
such that the mild solution $x\in \mathcal{C}(0,T;X)$ of system 
\eqref{system1} satisfies $x_{u}(T)=x_{d}$.	
\end{Definition}

To prove controllability, we make use of the following
assumptions $\left(A_1\right)$ and $\left(A_2\right)$:
\begin{enumerate}
\item[$\left(A_1\right)$]  
$\mathcal{T}(t)$ is compact for every $\displaystyle t>0$;

\item[$\left(A_2\right)$] 
The function $h:[0,T]\times \mathcal{C}(0, T ; X)\rightarrow X$ 
is continuous, and there exists a constant $\displaystyle \varsigma\in]0,T[$ 
and $\mathit{H},\mathit{H}_{1}>0$ such that $h\in D(\mathcal{A}^{\varsigma})$, 
and for any $z,y\in \mathcal{C}(0, T ; X)$, $t\in [0,T]$, the~function 
$\mathcal{A}^{\varsigma}h(\cdot,z)$ is strongly measurable and 
$\mathcal{A}^{\varsigma}h(t,\cdot)$ satisfies the Lipschitz condition
\begin{equation}
\label{eq8}
\displaystyle \left\vert \mathcal{A}^{\varsigma} h(t, z)
-\mathcal{A}^{\varsigma} h(t, y)\right\vert \leq \mathit{H}\|z-y\|
\end{equation}
and
\begin{equation}
\label{eq9}
\displaystyle \left\vert \mathcal{A}^{\varsigma} h(t, z)\right\vert 
\leq \mathit{H}_{1}\left(\|z\|+1\right).
\end{equation}
\end{enumerate}

\medskip

Let $\displaystyle H_{\nu}: L^{2}(0,T; U) \rightarrow X$ 
be the linear operator defined by 
$$
\displaystyle H_{\nu}u=\int_{0}^{T}(T-s)^{\nu-1} 
\mathit{K}_{\nu}(T-s) \mathcal{B}u(s)\mathrm{d} s.
$$ 
By construction, this operator is invertible. 
Indeed, because~$\displaystyle H_{\nu}$ takes values in the cokernel 
$L^{2}(0,T; U) \big/ {ker} H_{\nu}$, then it is injective. 
It is also surjective because 
$\displaystyle L^{2}(0,T; U) \big/ {ker} H_{\nu}\simeq Im H_{\nu}$ 
(see~\cite{ref1,ref5}). The~inverse operator 
$\displaystyle H_{\nu}^{-1}$ takes values in 
$\displaystyle L^{2}(0,T; U) \big/ \operatorname{ker} H_{\nu}$.
Thus, there exists a positive constant $\displaystyle M_{2}\geq 0$ 
such that
\begin{equation}
\label{eq10}
\displaystyle  \left\vert H_{\nu}^{-1}\right\vert_{\mathcal{L}\left(X, L^{2}(0,T; U) 
\big/\mathrm{ker}H_{\nu}\right)} \leq M_{2}.
\end{equation}

Let $\displaystyle r\geq 0$. Note that 
$\displaystyle B_{r}=\{x \in \mathcal{C}(0, T; X):\|x\| \leq r\}$  
is a bounded closed and convex subset in 
$\displaystyle \mathcal{C}(0, T; X)$.

\begin{Theorem}
\label{theorem 3.2}
If $(A_1)$ and $(A_2)$ are fulfilled, 
then the evolution system $\eqref{system1}$ is controllable 
in \linebreak  $\displaystyle [0,T]$ provided
\begin{equation}
\label{eq:Thm:1}
\begin{aligned}
\bigg[\vert \mathcal{A}^{-\varsigma}\vert+\frac{\mathbb{L}_{1-\varsigma}
\Gamma(1+\varsigma)}{\varsigma\Gamma(1+\nu\varsigma)}T^{\nu\varsigma}
+\frac{M\mathit{M}_{T}\mathit{M}_{1}}{\Gamma(1+\nu)}T^{\nu} 
\bigg(\vert \mathcal{A}^{- \varsigma}\vert 
+\displaystyle\frac{\mathbb{L}_{1-\varsigma}
\Gamma(1+\varsigma)}{\varsigma\Gamma(1+\nu\varsigma)}
T^{\nu\varsigma}\bigg)\bigg]\mathit{H}<1.
\end{aligned}
\end{equation}
\end{Theorem}

\begin{proof}
For any function $x$, we define the control
\begin{equation}
\label{eq12}
\begin{aligned}
u_{x}(t)& =H_{\nu}^{-1}\Big [x_{d}-\mathit{S}_{\nu}(t)\left[
x_{0}-h\left(0, x_{0}\right)\right]-h\left(T, x_{T}\right)\\
&\quad -\displaystyle\int_{0}^{T}(T-s)^{\nu-1} \mathcal{A} 
\mathit{K}_{\nu}(T-s) h\left(s, x_{s}\right) \mathrm{d} s\Big](t).
\end{aligned}
\end{equation}
We shall prove that $\mathit{G}:\mathcal{C}(0,T;X)
\rightarrow \mathcal{C}(0,T;X)$, defined by
\begin{equation}
\label{eqG}
\begin{aligned} 
(\mathit{G}x(t)) 
&=\mathit{S}_{\nu}(t)\left[x_{0}-h\left(0, x_{0}\right)\right]
+h\left(t, x_{t}\right)+\int_{0}^{t}(t-s)^{\nu-1} \mathcal{A} 
\mathit{K}_{\nu}(t-s) h\left(s, x_{s}\right) \mathrm{d} s\\ 
&+\int_{0}^{t}(t-s)^{\nu-1} \mathit{K}_{\nu}(t-s) \mathcal{B}
u_{x}(s)\mathrm{d} s,
\quad t \in[0, T], 
\end{aligned}
\end{equation}
has a fixed point $x$ for the control $u_{x}$ steering system 
$\eqref{system1}$ from $x_{0}$ to $x_{d}$ in time $T$. \linebreak  
From $\eqref{eq2}$, $\eqref{eq10}$, Lemma~\ref{YF} and (i) 
of Lemma~\ref{2.2}, we have
$$
\begin{aligned}
\left\vert \mathcal{B}u_{x}(t)\right\vert 
&\leq M\mathit{M}_{1}\bigg(\vert x_{d}\vert+\mathit{M}_{T}\Big[
\vert x_{0}\vert+\vert h(0,x_{0})\vert\Big]+\vert h(T,x_{T})\vert\\
&\quad +\displaystyle\int_{0}^{T}(T-s)^{\nu-1}\left\vert
\mathcal{A}^{1-\varsigma}\mathit{K}_{\nu}(T-s)
\mathcal{A}^{\varsigma}h(s,x_{s})\right\vert ds\bigg).
\end{aligned}
$$
In view of $\eqref{eq9}$ and (ii) of Lemma~\ref{2.2}, it follows that
$$
\begin{aligned}
\left\vert \mathcal{B}u_{x}(t)\right\vert 
&\leq M\mathit{M}_{1}\Bigg(\vert x_{d}\vert+\mathit{M}_{T}\Big[\vert x\vert
+\Big(r+1\Big)\mathit{H}_{1}\vert \mathcal{A}^{-\varsigma}\vert\Big]
+\Big(r+1\Big)\mathit{H}_{1}\vert \mathcal{A}^{-\varsigma}\vert\\
&\quad +\frac{\nu\mathbb{L}_{1-\varsigma}\Gamma(1+\varsigma)}{\Gamma(1+\nu\varsigma)}
\mathit{H}_{1}\Big(r+1\Big)\displaystyle\int_{0}^{T}(T-s)^{\nu\varsigma-1}ds\Bigg)\\
&\leq M\mathit{M}_{1}\Bigg(\vert x_{d}\vert+\mathit{M}_{T}\Big[\vert x\vert
+\Big(r+1\Big)\mathit{H}_{1}\vert \mathcal{A}^{-\varsigma}\vert \Big]
+\Big(r+1\Big)\mathit{H}_{1}\vert \mathcal{A}^{-\varsigma}\vert\\
&\quad +\displaystyle\frac{\mathbb{L}_{1-\varsigma}}{\Gamma(1+\varsigma)}
\mathit{H}_{1}\Big(r+1\Big)T^{\nu\varsigma}\Bigg).
\end{aligned}
$$
Let \vspace{-6pt}
\begin{multline*}
\mathcal{Y}=M\mathit{M}_{1}\Bigg(\vert x_{d}\vert
+\mathit{M}_{T}\left[\vert x\vert 
+\Big(r+1\Big)\mathit{H}_{1}\vert \mathcal{A}^{-\varsigma}\vert \right]\\
+\Big(r+1\Big)\mathit{H}_{1}\vert \mathcal{A}^{-\varsigma}\vert
+\frac{\mathbb{L}_{1-\varsigma}}{\Gamma(1+\varsigma)}
\mathit{H}_{1}\Big(r+1\Big)T^{\nu\varsigma}\Bigg).
\end{multline*}
It follows that
\begin{equation}
\label{eq13}
\left\vert \mathcal{B}u_{x}(t)\right\vert\leq \mathcal{Y}.
\end{equation}
In order to show that $\mathit{G}$ has a unique fixed point on $B_{r}$, 
we will proceed in two~steps.
 
\noindent Step I: $\mathit{G}x\in B_{r}$ whenever $x\in B_{r}$.
For any fixed $x\in B_{r}$ and $0\leq t\leq T$, we have\vspace{-6pt}
$$
\begin{aligned}
\left\vert(\mathit{G}x(t))\right\vert
\leq& \left\vert \mathit{S}_{\nu}(t)[x_{0}-h(0,x_{0})]\right\vert
+\vert h(t,x_{t})\vert +\int_{0}^{t}\left\vert (t-s)^{\nu-1}
\mathcal{A}\mathit{K}_{\nu}(t-s)h(s,x_{s})\right\vert~ds\\
&+\int_{0}^{t}(t-s)^{\nu-1}\vert \mathit{K}_{\nu}(t-s)\mathcal{B}u_{x}(s)\vert ds.
\end{aligned}
$$
From Lemma~\ref{YF}, \eqref{eq9}, and~(i) of Lemma~\ref{2.2}, it results that
$$
\begin{aligned}
\left\vert(\mathit{G}x(t))\right\vert
&\leq \mathit{M}_{T}\Big[r+\Big(r+1\Big)\mathit{H}_{1}\vert 
\mathcal{A}^{-\varsigma}\vert\Big]+\Big(r+1\Big)
\mathit{H}_{1}\vert \mathcal{A}^{-\varsigma}\vert\\
&\quad +\int_{0}^{t}(t-s)^{\nu -1}\Big\vert \mathcal{A}^{1-\varsigma}
\mathit{K}_{\nu}(t-s)\mathcal{A}^{\varsigma}h(s,x_{s})\Big\vert~ds\\
&\quad +\frac{\nu \mathit{M}_{T}}{\Gamma(1+\nu)}\int_{0}^{t}(t-s)^{\nu-1}
\vert \mathcal{B}u_{x}\vert ds.
\end{aligned}
$$
Now, by~using (ii) of Lemma~\ref{2.2}, we get 
$$
\begin{aligned}
\left\vert(\mathit{G}x(t))\right\vert
&\leq \mathit{M}_{T}[r+\mathit{H}\vert 
\mathcal{A}^{-\varsigma}\vert\Big(r+1\Big)]
+\mathit{H}\vert \mathcal{A}^{-\varsigma}\vert\Big(r+1\Big)\Big\vert\\
&\quad + \displaystyle\frac{\nu \mathbb{L}_{1-\varsigma}
\Gamma(1+\varsigma)}{\Gamma(1+\nu\varsigma)}
\mathit{H}\Big(r+1\Big)\int_{0}^{t}(t-s)^{\nu\varsigma-1}ds\\
&\quad +\displaystyle\frac{\nu \mathit{M}_{T}}{\Gamma(1+\nu)}
\int_{0}^{t}(t-s)^{\nu-1}\Big\vert \mathcal{B}u_{x}(s)\Big\vert ds.
\end{aligned}
$$
According to \eqref{eq13}, one has
$$
\begin{aligned}
\left\vert(\mathit{G}x(t))\right\vert
&\leq \mathit{M}_{T}\Bigg[r+\mathit{H}\vert 
\mathcal{A}^{-\varsigma}\vert\Big(r+1\Big)\Bigg]
+\mathit{H}\vert \mathcal{A}^{-\varsigma}\vert\Big(r+1\Big)\vert\\
&\quad + \frac{\nu \mathbb{L}_{1-\varsigma}\Gamma(1+\varsigma)}{\varsigma
\Gamma(1+\nu\varsigma)}\mathit{H}\Big(r+1\Big)T^{\nu\varsigma}
+\frac{\mathit{M}_{T}}{\Gamma(1+\nu)}\mathcal{Y} T^{\nu}.
\end{aligned}
$$	
By choosing
$$
\begin{aligned}
r&=\mathit{M}_{T}\Bigg[r+\Big(r+1\Big)\mathit{H}_{1}\vert
\mathcal{A}^{-\varsigma}\vert\Bigg]+\Big(r+1\Big)
\mathit{H}_{1}\vert \mathcal{A}^{-\varsigma}\vert\\
&\quad +\displaystyle \frac{\nu \mathbb{L}_{1-\varsigma}
\Gamma(1+\varsigma)}{\varsigma\Gamma(1+\nu\varsigma)}
\mathit{H}_{1}\Big(r+1\Big)T^{\nu\varsigma}
+\displaystyle\frac{\nu \mathit{M}_{T}}{\Gamma(1+\nu)}\mathcal{Y} T^{\nu},
\end{aligned}
$$
we get that $\mathit{G}x\in B_{r}$ whenever $x\in B_{r}$.

\noindent Step II: $\mathit{G}$ is a contraction on $B_{r}$.
For any $v,w\in B_{r}$ and $0\leq t\leq T$, 
in accordance with \eqref{eq12}, we obtain\vspace{-6pt}
$$
\begin{aligned}
\left\vert(\mathit{G} v)(t)-(\mathit{G} w)(t)\right\vert 
&\leq \bigg\vert h(t,v_{t})-h(t,w_{t})\bigg\vert\\
& +\displaystyle\int_{0}^{t}(t-s)^{\nu-1}\bigg\vert 
\mathcal{A}\mathit{r}_{\nu}(t-s)\bigg(h\big(s, 
v(s)\big)-h\big(s, w(s)\big)\bigg)\bigg\vert~ds \\
&+\displaystyle\int_{0}^{t}(t-s)^{\nu-1} \bigg\vert 
\mathit{r}_{\nu}(t-s) \mathcal{B}
H_{\nu}^{-1}\left[h(T,v_{T})-h(T,w_{T})
+\int_{0}^{T}(T-\tau)^{\nu-1}\right.\\
&\left.\times \mathcal{A} \mathit{K}_{\nu}(T
-\tau)\bigg(h\big(\tau, 
v(\tau)\big)-h\big(\tau, w(\tau)\big)\bigg) 
d\tau\right](s) \bigg\vert d s.
\end{aligned}
$$
Considering Lemma~\ref{2.2} and $(A_2)$, we get 
$$
\begin{aligned}
\left\vert(\mathit{G} v)(t)-(\mathit{G} w)(t)\right\vert
&\leq \mathit{H}\vert \mathcal{A}^{-\varsigma}\left\vert v
-w\right\vert+\displaystyle\frac{\nu \mathbb{L}_{1-\varsigma} 
\Gamma(1+\varsigma)}{\Gamma(1+\nu\varsigma)}\mathit{H}\left\vert 
v-w\right\vert\displaystyle\int_{0}^{t}(t-s)^{\nu\varsigma-1}ds\\
&+\displaystyle\frac{\nu M\mathit{M}_{T}\mathit{M}_{1}}{\Gamma(1+\nu)}
\displaystyle\int_{0}^{t}(t-s)^{\nu-1}\bigg[
\bigg\vert h(T,v_{T})-h(T,w_{T})\bigg\vert\\
&+\int_{0}^{t}(T-\tau)^{\nu-1}\bigg\vert \mathcal{A}^{1-\varsigma}\mathit{K}_{\nu}(t-\tau)
\mathcal{A}^{\varsigma}\bigg[h(\tau,v(\tau))-h(\tau,w(\tau))\bigg]\bigg\vert d\tau\bigg]ds.
\end{aligned}
$$
From \eqref{eq8}, we obtain that
$$
\begin{aligned}
\left\vert(\mathit{G} v)(t)-(\mathit{G} w)(t)\right\vert
&\leq \mathit{H}\vert \mathcal{A}^{-\varsigma}\left\vert v-w\right\vert
+\displaystyle\frac{\mathbb{L}_{1-\varsigma} \Gamma(1+\varsigma)}{\varsigma
\Gamma(1+\nu\varsigma)}\mathit{H}\left\vert v-w\right\vert T^{\nu\varsigma}\\
&\quad +\displaystyle\frac{\nu M\mathit{M}_{T}\mathit{M}_{1}}{\Gamma(1+\nu)}
\displaystyle\int_{0}^{T}(t-s)^{\nu-1}\bigg[\mathit{H}\vert 
\mathcal{A}^{-\varsigma}\left\vert v-w\right\vert\\
&\quad +\displaystyle\frac{\mathbb{L}_{1-\varsigma}\Gamma(1+\varsigma)}{\varsigma
\Gamma(1+\nu\varsigma)}\mathit{H}\left\vert v-w\right\vert T^{\nu\varsigma}\bigg]ds\\
&\leq \mathit{H}\vert \mathcal{A}^{-\varsigma}\left\vert v-w\right\vert
+\displaystyle\frac{\mathbb{L}_{1-\varsigma} \Gamma(1+\varsigma)}{\varsigma
\Gamma(1+\nu\varsigma)}\mathit{H}\left\vert v-w\right\vert T^{\nu\varsigma}\\
&\quad +\displaystyle\frac{M\mathit{M}_{T}\mathit{M}_{1}}{\Gamma(1+\nu)}T^{\nu}\bigg[
\vert \mathcal{A}^{-\varsigma}\vert
+\displaystyle\frac{\mathbb{L}_{1-\varsigma}\Gamma(1+\varsigma)}{\varsigma
\Gamma(1+\nu\varsigma)}T^{\nu\varsigma}\bigg]\mathit{H}\left\vert v-w\right\vert\\
&=\bigg[\vert \mathcal{A}^{-\varsigma}\vert+\frac{\mathbb{L}_{1-\varsigma}
\Gamma(1+\varsigma)}{\varsigma\Gamma(1+\nu\varsigma)}T^{\nu\varsigma}
+\frac{M\mathit{M}_{T}\mathit{M}_{1}}{\Gamma(1+\nu)}T^{\nu}\\
&\quad \bigg(\vert \mathcal{A}^{-\varsigma}\vert +\displaystyle
\frac{\mathbb{L}_{1-\varsigma}\Gamma(1+\varsigma)}{\varsigma \Gamma(1+\nu\varsigma)}
T^{\nu\varsigma}\bigg)\bigg]\mathit{H}\left\vert v-w\right\vert.
\end{aligned}
$$
From Theorem~\ref{theorem 3.2}, we have
$$
\bigg[\vert \mathcal{A}^{-\varsigma}\vert+\frac{\mathbb{L}_{1-\varsigma}
\Gamma(1+\varsigma)}{\varsigma\Gamma(1+\nu\varsigma)}T^{\nu\varsigma}
+\frac{M\mathit{M}_{T}\mathit{M}_{1}}{\Gamma(1+\nu)}
T^{\nu}\bigg(\vert \mathcal{A}^{-\varsigma}\vert
+\displaystyle\frac{\mathbb{L}_{1-\varsigma}\Gamma(1+\varsigma)}{\varsigma
\Gamma(1+\nu\varsigma)}T^{\nu\varsigma}\bigg)\bigg]\mathit{H}<1;
$$
it follows that
$$
\left\vert(\mathit{G} v)(t)-(\mathit{G} w)(t)\right\vert<\vert v-w\vert,
$$
that is, $\mathit{G}$ is a contraction on $B_{r}$. 
We conclude from the Banach fixed-point theorem that $\mathit{G}$ 
has a unique fixed point $x$ in $\mathcal{C}(0,T;X)$.
Then, by~injecting $u_{x}$ in \eqref{mild solution}, we have\vspace{-6pt}
$$
\begin{aligned}
x_{u_{x}}(T) &=\mathit{S}_{\nu}(T)\left[x_{0}-h\left(0, x_{0}\right)\right]
+h\left(T, x_{T}\right)+\int_{0}^{T}(T-s)^{\nu-1} \mathcal{A} 
\mathit{K}_{\nu}(T-s) h\left(s, x_{s}\right) \mathrm{d} s \\ 
&\quad +\int_{0}^{T}(T-s)^{\nu-1} \mathit{K}_{\nu}(T-s) \mathcal{B}u_{x}(s)\mathrm{d} s,\\
&=\mathit{S}_{\nu}(T)\left[x_{0}-h\left(0, x_{0}\right)\right]
+h\left(T, x_{T}\right)+\int_{0}^{T}(T-s)^{\nu-1} \mathcal{A} 
\mathit{K}_{\nu}(T-s) h\left(s, x_{s}\right) \mathrm{d} s \\ 
&\quad +H_{\nu}H_{\nu}^{-1}\Big [x_{d}-\mathit{S}_{\nu}(T)\left[
x_{0}-h\left(0, x_{0}\right)\right]-h\left(T, x_{T}\right)\\
&\quad -\displaystyle\int_{0}^{T}(T-s)^{\nu-1} \mathcal{A} 
\mathit{K}_{\nu}(T-s) h\left(s, x_{s}\right) \mathrm{d} s\Big]\\
&=x_{d}
\end{aligned}
$$
and system \eqref{system1} is exactly controllable, 
which completes the proof.
\end{proof}

We have shown, under~assumptions $(A_1)$ and $(A_2)$, 
and with the help of Schauder's fixed-point theorem, 
that the neutral system \eqref{system1} is controllable
when condition \eqref{eq:Thm:1} holds. 
It would be interesting to clarify if the obtained control 
is unique in the sense that any control 
that allows reaching the state $x_{d}$ is such that the associated 
state $x$ is a fixed point of the operator $G$.
This uniqueness question is relevant but remains open.


\section{Optimal~Control}
\label{sec4}

Now, we consider the problem of steering system \eqref{system1} from the state $x_{0}$ 
to a target state $x_{d}$ in time $T$ with minimum energy. 
We prove the existence of solution to such an optimal control problem
when the set of admissible controls is closed and~convex.

Let $\mathcal{U}_{ad}$ be the nonempty set of admissible controls defined by 
$$
\mathcal{U}_{ad}=\left\{u \in L^{2}\left(0, T ; U\right): x_{u}(T)=x_{d}\right\}.
$$ 
\noindent We shall prove that $\mathcal{U}_{ad}$ is closed. For~that, 
let us consider a sequence $u_{n}$ in $\mathcal{U}_{ad}$ such that 
$u_{n}\rightarrow u$ strongly in $L^{2}(0,T;U)$, so
$$
\begin{aligned}
x_{u_{n}}(T)&=\mathit{S}_{\nu}(T)\left[x_{0}-h\left(0, x_{0}\right)\right]
+h\left(T, x_{T}\right)+\int_{0}^{T}(T-s)^{\nu-1} \mathcal{A} 
\mathit{K}_{\nu}(T-s) h\left(s, x_{s}\right) \mathrm{d} s \\ 
&\quad +\int_{0}^{T}(T-s)^{\nu-1} \mathit{K}_{\nu}(T-s) 
\mathcal{B}u_{n}(s)\mathrm{d} s.
\end{aligned}
$$
Put
$$
\mathcal{Q}u=\displaystyle\int_{0}^{T}(T-s)^{\nu-1} \mathcal{A} 
\mathit{K}_{\nu}(T-s) h\left(s, x_{s}\right) \mathrm{d} s
+\displaystyle\int_{0}^{T}(T-s)^{\nu-1} \mathit{K}_{\nu}(T-s) \mathcal{B}u_{n}(s)\mathrm{d} s.
$$
Since $\mathcal{Q}u$ is continuous, then $\mathcal{Q}u_{n}\rightarrow\mathcal{Q}u$ 
strongly in $X$. We also have that $h:[0,T]\times \mathcal{C}(0, T ; X)\rightarrow X$ 
is continuous; then $x_{u_{n}}(T)\rightarrow x_{u}(T)$ in $X$,
but $x_{u_{n}}(T)\in\{x_{d}\}$, which is closed. Therefore, $x_{u}(T)\in\{x_{d}\}$, 
which means that $u\in\mathcal{U}_{ad}$. Hence, $\mathcal{U}_{ad}$ is closed.

For a desired state $x_{d}$, our optimal control problem 
consists of finding within $\mathcal{U}_{ad}$ a control minimizing 
the functional
$$
J(u)=\frac{\varsigma}{2} \int_{0}^{T}\left\vert x_{u}(t)-x_{d}\right\vert_{X}^{2} d t
+\frac{\varepsilon}{2} \int_{0}^{T}\vert u(t)\vert_{U}^{2} d t,
$$
where $x_{u}(\cdot)$ is the mild solution of system \eqref{system1} 
associated with $u$. The~parameters $\varepsilon$ and $\varsigma$ 
are non-negative constants. Precisely, our optimal control problem is:
\begin{equation}
\label{eq14}
\left\{\begin{array}{l}
\displaystyle\inf_{u \in \mathcal{U}_{a d}} J(u),\\
\text{s.t. } \eqref{system1}.
\end{array}\right.
\end{equation}

The following result gives a necessary condition for the existence
of an optimal control to our minimum energy~problem.

\begin{Theorem}
\label{res2:Thm2}
Let $\mathcal{U}_{ad}$ be closed and convex.
If $1-\mathit{H}\left\vert \mathcal{A}^{-\varsigma}\right\vert>0$,
then there exists a $u^{\star}\in \mathcal{U}_{ad}$ solution 
to the optimal control problem \eqref{eq14}.
\end{Theorem}

\begin{proof}
Let $\displaystyle\left\vert u_{p}\displaystyle\right\vert^{2}
\leq\displaystyle\frac{2}{\displaystyle\varepsilon}\displaystyle J(u_{p})$ 
with $\displaystyle(u_{p})_{p\in \mathbb{N}}$ bounded.
Then there exists a subsequence, still denoted $(u_{p})_{p\in \mathbb{N}}$, 
that converges weakly to a limit $u^{\star}$. If~$\mathcal{U}_{ad}$ is
closed and convex, then $\mathcal{U}_{ad}$ is closed for the weak topology, 
which implies that $u^{\star}\in \mathcal{U}_{ad}$.
Let $x_{p}$ be the unique solution of system \eqref{system1} associated 
with $u_{p}$, and let $x^{\star}$ be the unique solution of system \eqref{system1} 
associated with $u^{\star}$. Then,
\begin{equation}
\begin{aligned}
\left\vert x_{p}(t)-x^{*}(t)\right\vert 
&\leq \left\vert h\left(t, x_{p}(t)\right)
-h\left(t, x^{\star}(t)\right)\right\vert\\
&\quad +\left\vert\int_{0}^{t}(t-s)^{\nu-1} \mathcal{A} 
\mathit{K}_{\nu}(t-s) [h\left(s, x_{p}(s)\right)
-h\left(s, x^{\star}(s)\right)] \mathrm{d} s\right\vert\\ 
&\quad +\left\vert\int_{0}^{t}(t-s)^{\nu-1} \mathit{K}_{\nu}(t-s) 
\mathcal{B}[u_{p}(s)-u^{\star}(s)]\mathrm{d} s\right\vert\\
&\leq \mathit{H}\left\vert \mathcal{A}^{-\varsigma}
\right\vert\left\vert x_{p}(t)-x^{\star}(t)\right\vert\\
&\quad +\int_{0}^{t}(t-s)^{\nu-1} \left\vert \mathcal{A}^{1-\varsigma} 
\mathit{K}_{\nu}(t-s) [\mathcal{A}^{\varsigma}
h\left(s, x_{p}(s)\right)-\mathcal{A}^{\varsigma}
h\left(s, x^{\star}(s)\right)]\right\vert \mathrm{d} s\\
&\quad +\left\vert\int_{0}^{t}(t-s)^{\nu-1} \mathit{K}_{\nu}(t-s) 
\mathcal{B}[u_{p}(s)-u^{\star}(s)]\mathrm{d} s\right\vert,
\quad t \in[0, T]. 
\end{aligned}
\end{equation}
This leads us to
\begin{equation}
\begin{aligned}
\left(1-\mathit{H}\left\vert \mathcal{A}^{-\varsigma}
\right\vert\right)\left\vert x_{p}(t)-x^{*}(t)\right\vert 
&\leq\frac{\nu \Gamma(1+\varsigma)}{\Gamma(1+\nu \varsigma)}
\mathbb{L}_{1-\varsigma}\int_{0}^{t}(t-s)^{\nu \varsigma-1}
\mathit{H}\left\vert x_{p}(t)-x^{\star}(t)\right\vert\mathrm{d}s\\
&\quad +\left\vert\int_{0}^{t}(t-s)^{\nu-1} \mathit{K}_{\nu}(t-s) 
\mathcal{B}[u_{p}(s)-u^{\star}(s)]\mathrm{d} s\right\vert,
\end{aligned}
\end{equation}
$t \in[0, T]$. Set $\mathcal{K}^{\prime}=\displaystyle\frac{1}{1
-\mathit{H}\left\vert \mathcal{A}^{-\varsigma}\right\vert}$. Then,
\begin{equation}
\begin{aligned}
\left\vert x_{p}(t)-x^{*}(t)\right\vert 
&\leq \mathcal{K}^{\prime}\frac{\nu \Gamma(1+\varsigma)}{
\Gamma(1+\nu \varsigma)}\mathbb{L}_{1-\varsigma}
\int_{0}^{t}(t-s)^{\nu \varsigma-1}
\mathit{H}\left\vert x_{p}(t)-x^{\star}(t)\right\vert\mathrm{d}s\\
&\quad +\mathcal{K}^{\prime}\left\vert\int_{0}^{t}(t-s)^{\nu-1} 
\mathit{K}_{\nu}(t-s) \mathcal{B}[u_{p}(s)-u^{\star}(s)]
\mathrm{d} s\right\vert,\quad t \in[0, T].
\end{aligned}
\end{equation}
Using the Gronwall lemma, we obtain that
\begin{equation}
\begin{aligned}
\left\vert x_{p}(t)-x^{*}(t)\right\vert 
&\leq \mathcal{K}^{\prime}\left\vert\int_{0}^{t}(t-s)^{\nu-1} 
\mathit{K}_{\nu}(t-s) \mathcal{B}[u_{p}(s)-u^{\star}(s)]\mathrm{d}s
\right\vert\\
&\quad \exp\left( \mathcal{K}^{\prime}\frac{\nu 
\Gamma(1+\varsigma)}{\Gamma(1+\nu \varsigma)}
\mathbb{L}_{1-\varsigma}\mathit{H}
\int_{0}^{t}(t-s)^{\nu \varsigma-1}\mathrm{d}s\right)\\
&\leq \mathcal{K}^{\prime}\left\vert\int_{0}^{t}(t-s)^{\nu-1} 
\mathit{K}_{\nu}(t-s) \mathcal{B}[u_{p}(s)-u^{\star}(s)]
\mathrm{d} s\right\vert\\
&\quad \exp\left(\mathcal{K}^{\prime}\frac{ \Gamma(1+\varsigma)}{
\varsigma\Gamma(1+\nu \varsigma)}\mathbb{L}_{1-\varsigma}
\mathit{H} T^{\nu\varsigma}\right).
\end{aligned}
\end{equation}
Now, by~the weak convergence, $u_{p} \rightharpoonup u^{*}$ 
in ${L^{2}(0,T,U)}$, and~from Lemma~\ref{YF}, we obtain that
\begin{multline}
\left\vert \int_{0}^{t}(t-s)^{\nu-1} \mathit{K}_{\nu}(t-s) 
\mathcal{B}[u_{p}(s)-u^{\star}(s)]\mathrm{d} s\right\vert\\
\leq\frac{\nu \mathit{M}_{T}\mathit{M}_{1}}{\Gamma(1+\nu)}
\int_{0}^{t}(t-s)^{\nu-1}
\left\vert u_{p}(s)-u^{\star}(s)\right\vert_{L^{2}(0,T,U)}\mathrm{d} s,
\end{multline}
from which 
$x_{p}\rightarrow x^{\star}$ strongly in $L^{2}(0,T;X)$. Hence,
$$
\lim_{n\rightarrow\infty}\int_{0}^{T}\left\vert 
x_{p}(t)-x_{d}\right\vert_{X}^{2} d t
=\int_{0}^{T}\left\vert x(t)-x_{d}\right\vert_{X}^{2} d t.
$$
Using the lower semi-continuity of norms, 
the weak convergence of $(u_{p})_{n}$ gives
$$
\left\vert u^{\star}\right\vert
\leq \displaystyle \lim_{n\rightarrow\infty}
\inf\left\vert u_{p}\right\vert.
$$
Therefore, $J(u^{\star})\leq \displaystyle\lim_{n\rightarrow\infty}\inf J(u_{p})$, 
leading to $J(u^{\star})=\displaystyle\inf_{u\in \mathcal{U}_{ad}}J(u_{p})$,
which establishes the optimality of $u^{\star}$.
\end{proof}

We have just proved the existence of an optimal control for a closed convex set 
of admissible controls. In~Section~\ref{sec5}, our main results are 
illustrated with the help of an example.


\section{An~Application}
\label{sec5}

In this section we illustrate the results given
by our Theorems~\ref{theorem 3.2} and \ref{res2:Thm2}.

Let $\displaystyle X= L^{2}((0,1);\mathbb{R})$ and consider 
the fractional differential system
\begin{equation}
\label{ex}
\left\{
\begin{array}{clclcl}
\displaystyle { }^{C} D_{t}^{1/2}\Big(y(t,z) -h(t,y_{t})\Big) 
= \Delta y(t,z)+\mathcal{B}u(t,z), & t\in [0,1],\\
y(t,0) = y(t,1)=0, & t\in[0,1], 
\end{array}
\right.
\end{equation}
where the order $\nu$ of the fractional derivative 
is equal to $\displaystyle\frac{1}{2}$, and~the function 
$\displaystyle h:[0,1]\times\mathcal{C}\rightarrow X$ is given by
\begin{equation}
\label{funch}
h(t,y_{t})(x)=\int_{0}^{1}\mathcal{F}(x,z)u_{t}(v,z)dz,
\end{equation}
where $\mathcal{F}$ is assumed to satisfy the following conditions:
\begin{enumerate}
\item[(a)] The function $\displaystyle \mathcal{F}(x,z)$, $x,z\in[0,1]$, 
is measurable and
$$
\displaystyle \int_{0}^{1}\int_{0}^{1}\mathcal{F}^{2}(x,z)dz<\infty;
$$

\item[(b)] The function $\displaystyle \partial{x}\mathcal{F}(x,z)$ is measurable, 
$\displaystyle  \mathcal{F}(0,z)=\mathcal{F}(1,z)=0$, and~
$$
\displaystyle \Bigg(\int_{0}^{1}\int_{0}^{1}\big(\partial{x}
\mathcal{F}(x,z)\big)^{2}dzdx\Bigg)^{1/2}<\infty.
$$
\end{enumerate}

Let $\mathcal{A}:D(\mathcal{A})\subseteq X\rightarrow X$ be defined by 
$\displaystyle \mathcal{A}x=-x^{\prime\prime}$ with the domain
$$
\displaystyle D(\mathcal{A})=\left\{x(\cdot)\in X : 
x, x^{\prime} ~\mbox{ absolutely continuous }, 
x^{\prime\prime}\in X, x(0)=x(1)=0\right\}.
$$

We begin by proving that the assumption $(A_1)$ holds.
Indeed, operator $\displaystyle \mathcal{A}$ is self-adjoint, 
with a compact resolvent, and~generating an analytic compact semi-group 
$\mathcal{T}(t)$. Furthermore, the~eigenvalues of $\mathcal{A}$ 
are $\displaystyle \Lambda_{p}=p^{2}\pi^{2}$, $p\in \mathbb{N}$, 
with corresponding normalized eigenvectors 
$\displaystyle e_{p}(z)=\sqrt{\frac{2}{\pi}}\sin(p\pi z)$, 
$\{e_{i}\}_{i=1}^{\infty}$ forming an orthonormal basis of $\displaystyle X$.
Then, \vspace{-6pt}
$$
\mathcal{A}x=-\sum_{p=1}^{p=\infty}\Lambda_{p}(x,e_{p})e_{p},
\quad x\in D(\mathcal{A}),
$$
and 
$$
\mathcal{T}(t)x(s)=\sum_{i=1}^{i=\infty}
\exp(\Lambda_{i}t)(x,e_{i})e_{i}(s),
\quad x\in X.
$$

Note that $\displaystyle \mathcal{T}(\cdot)$ is a uniformly stable semi-group 
and $\displaystyle \|\mathcal{T}(t)\|_{L^{2}[0,1]}\leq \exp(-t)$. 
The following properties~hold:
\begin{enumerate}
\item[(i)] $\displaystyle \mathcal{A}^{-\frac{1}{2}}x
=\sum_{p=1}^{\infty}\frac{1}{p}(x,e_{p})e_{p}$;

\item[(ii)] The operator $\mathcal{A}^{\frac{1}{2}}$ is given by\vspace{-6pt}
$$
\displaystyle \mathcal{A}^{\frac{1}{2}}x=\sum_{p=1}^{\infty}p(x,e_{p})e_{p}
$$
and $\displaystyle D(\mathcal{A}^{\frac{1}{2}})
=\left\{x(\cdot)\in X,\sum_{p=1}^{\infty}p(x,e_{p})e_{p}\in X\right\}$.
\end{enumerate}
Clearly, \eqref{M_{T}}, \eqref{C}, and~
$(A_1)$ are~satisfied.

Under our assumptions (a) and (b) on $\mathcal{F}$,
\eqref{eq8} and \eqref{eq9} are also satisfied,
and assumption $(A_2)$ also holds.

Let $U$ be a reflexive Banach space. We consider the
control operator $\displaystyle \mathcal{B}:U\rightarrow X$ defined by\vspace{-6pt}
$$
\mathcal{B}u=\sum_{p=1}^{p=\infty}\Lambda_{p}({\bar{u}},e_{p})e_{p},
$$
where 
$$
\bar{u}
= \left\{
\begin{array}{clclcl}
u_{p}, \quad &p=1,2,\ldots N,\\
0, \quad & p=N+1,N+2,\ldots
\end{array}
\right.
$$
We see that $\mathcal{B}$ is a bounded continuous operator 
with $\displaystyle \mathit{M}_{1}=N\Lambda_{N}$.
For $\displaystyle N\in \mathbb{N}$ 
and $\displaystyle H_{1/2}:L^{2}([0,1],U)\rightarrow X$ given by\vspace{-6pt}
$$
\displaystyle H_{1/2}u=\int_{0}^{1}(1-s)^{1/2}P_{1/2}(1-s)\mathcal{B}u(s)ds,
$$
we have\vspace{-6pt}
$$
\fontsize{9.5}{9.5}\selectfont\begin{aligned}
\displaystyle H_{1/2}u
&=\displaystyle\int_{0}^{1}(1-s)^{1/2}\frac{1}{2}
\displaystyle\int_{0}^{\infty}\Theta\phi_{1/2}(\Theta)
T((1-s)^{1/2}\Theta)\mathcal{B}u(s)d\Theta~ds\\
&=\displaystyle\int_{0}^{1}(1-s)^{1/2}\frac{1}{2}
\displaystyle\int_{0}^{\infty}\Theta\phi_{1/2}(\Theta)
\sum_{i=1}^{i=\infty}\exp(\Lambda_{i}(1-s)^{1/2}\Theta)
(\mathcal{B}u,e_{i})e_{i}(s)d\Theta~ds\\
&=\int_{0}^{1}(1-s)^{1/2}\sum_{i=1}^{\infty}\int_{0}^{\infty}
\frac{1}{2}\Theta\phi_{1/2}(\Theta)\sum_{j=0}^{\infty}
\frac{\Lambda_{i}(1-s)^{1/2}\Theta)^{j}}{j!}(u,e_{i})e_{i}(s)d\Theta~ds\\
&=\int_{0}^{1}(1-s)^{1/2}\sum_{i=1}^{\infty}\sum_{j=0}^{\infty}
\frac{(\Lambda_{i}(1-s)^{1/2})^{j}}{\Gamma(1/2+\frac{1}{2}j)}(u,e_{i})e_{i}(s)ds\\
&=\sum_{i=1}^{\infty}\sum_{j=0}^{\infty}\int_{0}^{1}
\frac{\Lambda_{i}^{j}}{\Gamma(\frac{1}{2}
+\frac{1}{2}j)}(1-s)^{\frac{1+j}{2}}(u,e_{i})e_{i}(s)\\
&=\sum_{i=1}^{\infty}\sum_{j=0}^{\infty}
\frac{2\Lambda_{i}^{j}}{\Gamma(\frac{1}{2}
+\frac{1}{2}j)(3+j)}(u,e_{i})e_{i}(s).
\end{aligned}
$$
Applying Theorem~\ref{theorem 3.2}, we deduce that the 
fractional differential system \eqref{ex} is controllable. 
Moreover, for~function $h$ defined as in \eqref{funch}
with the Lipshitz constant
$\mathit{H}<\displaystyle\frac{1}{\left\vert \mathcal{A}^{-\frac{1}{2}}\right\vert}$,
we conclude from Theorem~\ref{res2:Thm2} 
that there exists a control steering the system, in~one unit of time,
from a given initial state to a given terminal state with minimum~energy. 


\section{Conclusions}
\label{sec6}

Using the Banach fixed-point theorem, we
have obtained a set of sufficient conditions for the controllability 
of a class of fractional neutral evolution equations involving 
the Caputo fractional derivative of order $\alpha\in ]0,1[$ 
(cf. Theorem~\ref{theorem 3.2}). The result is proved in two major steps: 
(i) in the first step, we proved that the operator $G$ defined by \eqref{eqG} 
is an element of the bounded closed and convex subset $B_{r}$,
(ii) while in the second, we proved that $G$ is a contraction 
on the same subset $B_{r}$. Moreover, we formulated 
a minimum energy optimal control problem and proved
conditions assuring the existence of a solution 
for the optimal control problem
$\displaystyle\inf_{u \in \mathcal{U}_{a d}} J(u)$
subject to \eqref{system1} (cf. Theorem~\ref{res2:Thm2}).
An example was given illustrating the two main results.

Our work can be extended in several directions:
(i) to a case of enlarged controllability using different fractional 
derivatives; (ii) by developing methods to determine
the control predicted by our existence theorem, e.g.,
by using RHUM and penalization approaches~\cite{37,38,Karite21}; 
(iii) or by giving applications of neutral systems to epidemiological 
problems~\cite{Khan22,Zarin21}. Many other questions remain open, 
as~is the case of regional controllability and regional discrete 
controllability for problems of the type considered here. 
A strong motivation behind~the investigation of neutral evolution systems, 
such as \eqref{system1} considered here, comes from physics, since they describe 
well various physical phenomena as~fractional diffusion equations. 
However, neutral systems are difficult to study, since such control systems 
contain time-delays not only in the state but also in the velocity variables, 
which make them intrinsically more complicated. The~limitations 
of the method we proposed here is that we are not able
to provide conditions under which the optimal control 
is unique. Additionally, we do not have an explicit form for it.


\vspace{6pt} 


\authorcontributions{Conceptualization, A.A.; 
methodology, Z.E.-c., A.A., T.K. and D.F.M.T.; 
validation, Z.E.-c., A.A., T.K. and D.F.M.T.; 
formal analysis, Z.E.-c., A.A., T.K. and D.F.M.T.; 
investigation, Z.E.-c., A.A., T.K. and D.F.M.T.; 
writing---original draft preparation, Z.E.-c., A.A., T.K. and D.F.M.T.; 
writing---review and editing, Z.E.-c., A.A., T.K. and D.F.M.T. 
All authors have read and agreed to the published version of the~manuscript.}

\funding{T.K. and D.F.M.T. were partially funded by FCT, project UIDB/04106/2020 (CIDMA).}

\institutionalreview{Not applicable.}

\informedconsent{Not applicable.}

\dataavailability{Not applicable.} 

\acknowledgments{This research is part of Ech-chaffani's Ph.D., 
which is being carried out at Sidi Mohamed Ben Abdellah, Fez,
under the scientific supervision of Aberqi. It was essentially
finished during a one-month visit of Karite to the Department 
of Mathematics of University of Aveiro, Portugal, April and May 2022.
The authors are very grateful to two anonymous referees
for many suggestions and invaluable~comments.}

\conflictsofinterest{The authors declare no conflict of interest.
The funders had no role in the design of the study; in the collection, 
analyses, or~interpretation of data; in the writing of the manuscript; 
or in the decision to publish the~results.} 
\clearpage 


\end{paracol}

\reftitle{References}


\end{document}